\documentclass[11pt,fullpage]{article}

\usepackage{graphicx}
\usepackage{latexsym}
\usepackage{color}
\usepackage{multirow}
\usepackage{rotating}
\usepackage{amsmath}
\usepackage{graphicx}
\usepackage{amssymb}
\usepackage{amsthm}
\usepackage{amsfonts}

\renewcommand{\baselinestretch} {1.3}

\makeatletter 
\def\singlespace{\def\baselinestretch{1}\@normalsize}

\newtheorem{theorem}{Theorem}

\newtheorem{proposition}{Proposition}

\newcommand{\Jo}{{\cal J}_0}
\newcommand{\Jc}{{\cal J}_0^c}
\newcommand{\Fj}{{\cal F}_{n_j}(s_j,R_j)}
\newcommand{\Tj}{\Theta_{n_j}(s_j,R_j)}

\long\def\ignore#1{}

\newcommand{\be}{\begin{equation}}
\newcommand{\ee}{\end{equation}}
\newcommand{\beqn}{\begin{eqnarray}}
\newcommand{\eeqn}{\end{eqnarray}}

\begin{document}


\title{\bf Sparse additive regression on a regular lattice}
\author{{\bf Felix Abramovich \hspace{2.0cm} Tal Lahav} \\
\hspace{1.1cm}felix@post.tau.ac.il \hspace{1.6cm} lahav.t@gmail.com\\
\hspace{0.9cm}Department of Statistics and Operations Research \\
Tel Aviv University \\
Tel Aviv 69978 \\
Israel
}

\date{}

\maketitle

\begin{abstract}
We consider estimation in a sparse additive regression model with the design
points on a regular lattice. We establish
the minimax convergence rates over Sobolev classes and propose a
Fourier-based rate-optimal estimator which is adaptive to the unknown
sparsity and smoothness of the response function.
The estimator is derived within Bayesian formalism but
can be naturally viewed as a penalized
maximum likelihood estimator with the complexity penalties on the number of nonzero
univariate additive components of the response and on the numbers of the nonzero
coefficients of their Fourer expansions. We compare it with several
existing counterparts and perform a short simulation study to demonstrate its
performance.
\end{abstract}


\bigskip

\section{Introduction} \label{sec:intr}
Consider a general nonparametric $d$-dimensional regression model,
where the design points are located on a
regular lattice of size $n_1 \times \ldots \times  n_d$ on $[0,1]^d$:
\be y(i_1/n_1,\ldots,i_d/n_d)=f(i_1/n_1,\ldots,i_d/n_d)+\epsilon(i_1/n_1,...,i_d/n_d),
\quad i_j=0,\ldots,n_j-1;\;j=1,\ldots d \label{eq:model0}
\ee
$\epsilon(i_1/n_1,...,i_d/n_d) \sim {\cal N}(0,\sigma^2)$ and are independent,
and the unknown response function $f:\mathbb{R}^d \rightarrow \mathbb{R}$ is
assumed to belong to a class of functions of certain smoothness.
Let $N=\prod_{j=1}^d n_j$ be the overall number of observations in the model
(\ref{eq:model0}).

In particular, a regular grid can be useful for design of experiments when
one has some prior belief on the relative relevance of
predictors. Thus, he can use a finer grid (larger $n_j$) for more
important variables and a coarse grid (smaller $n_j$) otherwise.

When $d$ is large, estimation of $f$ in (\ref{eq:model0}) suffers
severely from ``curse of dimensionality'' problem. A typical remedy
is to impose some addition structural constraints on $f$. One of the
common approaches is to consider the class of {\em additive} models
(Hastie \& Tibshirani, 1990), where the unknown $f$ can be
decomposed in a sum of $d$ univariate functions:
$f(x_1,\ldots,x_d)=\sum_{j=1}^d f_j(x_j)$. The original model
(\ref{eq:model0}) becomes then
\be
y(i_1/n_1,\ldots,i_d/n_d)=a_0+\sum_{j=1}^d
f_j(i_j/n_j)+\epsilon(i_1/n_1,...,i_d/n_d), \quad
i_j=0,\ldots,n_j-1;\;j=1,\ldots d. \label{eq:admodel}
\ee
To make the model (\ref{eq:admodel}) identifiable, we impose $\sum_{i=0}^{n_j-1}
f_j(i/n_j)=0$ for all $j=1,\ldots,d$. The goal is to estimate the
unknown global mean $a_0$ and the functions $f_j$'s.

Additive models have become a standard tool in multivariate
nonparametric regression and can be efficiently fitted by the
backfitting algorithm of Friedman \& Stuetzle (1981). However, in a
variety of modern high-dimensional statistical setups the number of
predictors $d$ may be still large relatively to the amount of
observed data. A key extra assumption then is {\em sparsity}, where
it is assumed that only a small fraction of $f_j$ in
(\ref{eq:admodel}) has a truly relevant impact on the response while
other $f_j=0$. Let $\Jo$ and $\Jc$ be the (unknown) subsets of
indices corresponding respectively to the zero and nonzero $f_j$.
The {\em sparse additive model} is
\be
y(i_1/n_1,\ldots,i_d/n_d)=a_0+\sum_{j \in \Jc}
f_j(i_j/n_j)+\epsilon(i_1/n_1,...,i_d/n_d), \quad
i_j=0,\ldots,n_j-1;\;j=1,\ldots d \label{eq:spad}
\ee
and $\sum_{i=0}^{n_j-1} f_j(i/n_j)=0,\;j \in \Jc$.

Expand each $f_j,\;j \in \Jc$ in the orthogonal discrete Fourier series assuming
for simplicity of exposition that all $n_j$ are odd:
$$
f_j(i/n_j)=\sum_{k=-(n_j-1)/2}^{(n_j-1)/2} c_{kj} e^{-\frac{2\pi I k i}{n_j}},
$$
where $I=\sqrt{-1}$ and discrete Fourier coefficients
\be
\label{eq:disfourier}
c_{kj}=\frac{1}{n_j}\sum_{i=0}^{n_j-1} f_j(i/n_j) e^{\frac{2\pi I k i}{n_j}}.
\ee
The identifiability condition $\sum_{i=0}^{n_j-1}f_j(i/n_j)=0$
implies $c_{0j}=0$.

One should make some assumptions on regularity properties of $f_j$.
We assume that the vector of discrete Fourier coefficients $c_j$ of
$f_j$ in (\ref{eq:disfourier}) belongs to a Sobolev ellipsoid
$\Tj=\{c_j:\sum_{k=-(n_j-1)/2}^{(n_j-1)/2} |c_{kj}|^2 |k|^{2s_j} \leq
R^2_j;\;c_{0j}=0\}$, where $s_j>1/2$ and $R_j<C_R$ for some constant
$C_R>0$, and denote the corresponding class of functions $f_j$ by
$\Fj$. The class $\Fj$ is a discrete analog of a Sobolev ball of
functions of smoothness $s_j$ with a radius $R_j$ (see, e.g.,
Korostelev \& Korosteleva, 2011, Section 10.5).

We establish the minimax rates of estimating $f$ in (\ref{eq:spad}),
where $f_j \in \Fj$. The corresponding rates for the case of $N$ distinct
points for each predictor $x_j$ were derived
in Raskutti, Wainwright \& Yu (2012). However, we
consider a design on the regular lattice, where there are $N/n_j$
repeated observations at each of $n_j$ grid points for every $x_j$.
It turns out that this difference affects the resulting minimax rates. 

In particular, we show that the average mean squared error $AMSE(\hat{f}_j,f_j)=\frac{1}{n_j}E||\hat{f}_j-f_j||^2_{n_j}$
for estimating a single univariate function
$f_j \in \Fj$ in the model (\ref{eq:spad}) at the design points,
where a general notation $||\cdot||_n$ is used for Euclidean norm in $\mathbb{R}^n$,
is of the order
\be \label{eq:amsej}
\min\left(N^{-\frac{2s_j}{2s_j+1}},\frac{n_j}{N}\right).
\ee
For sufficiently smooth $f_j$ with $2s_j+1 \geq \ln N/\ln n_j$, the rate in
(\ref{eq:amsej}) is
the standard minimax rate $N^{-2s_j/(2s_j+1)}$ for nonparametric estimation
of a univariate function from $\Fj$ (see, e.g., Korostelev \& Korosteleva, 2011,
Section 10.5),
but for $2s_j+1<\ln N/\ln n_j$ it corresponds to the parametric
rate of estimating $f_j$ at each grid point $i/n_j$ by simple averaging over the
corresponding $N/n_j$ replications. To understand this phenomenon recall that
in a standard nonparametric regression setup smoothing (local averaging over
neighbour points) is necessary to reduce
the variance. Although it introduces bias, the effect of the latter is
negligible
under smoothness assumptions on an unknown response function, while the benefits
of variance reduction are essential. As we have mentioned above, in
the considered case
there are $N/n_j$ repeated observations at each grid point $i/n_j$ and the
variance can already be reduced by their averaging without causing any bias.
On the other hand,
the grid might be too coarse to use neighbour points in smoothing since
the resulting bias becomes dominating in the bias-variance tradeoff for
nonsmooth $f_j$, where $2s_j+1<\ln N/\ln n_j$.

In particular, when all $n_j=N^{1/d}$ are equal, the minimax
$AMSE(\hat{f}_j,f_j)$ in (\ref{eq:amsej}) is of the order $N^{-r_j}$, where
$r_j=\max\left(\frac{2s_j}{2s_j+1},1-\frac{1}{d}\right)$ and the parametric
rate of averaging occurs when $2s_j+1<d$.

Furthermore, we prove that the overall minimax
$AMSE(\hat{f},f)=\frac{1}{N}E||\hat{f}-f||^2_N$ for the sparse additive
models with $d_0=|\Jc|$ nonzero $f_j$ is of the order
\be \label{eq:amse}
\max\left(\sum_{j \in \Jc} \min\left(N^{-\frac{2s_j}{2s_j+1}},\frac{n_j}{N}\right), \frac{d_0 \ln(d/d_0)}{N}\right).
\ee
The term $\sum_{j \in \Jc} \min\left(N^{-\frac{2s_j}{2s_j+1}},\frac{n_j}{N}\right)$ in (\ref{eq:amse})
is associated with the minimax rates of estimating $d_0$ nonzero univariate
functions in $\Fj,\;j \in \Jc$, while $\frac{d_0 \ln(d/d_0)}{N}$ corresponds
to the error of selecting a subset of $d_0$ nonzero elements out of $d$ and
appears in various related model selection setups (e.g., Abramovich \& Grinshtein, 2010, 2013;
Raskutti, Wainwright \& Yu, 2011, 2012; Rigollet \& Tsybakov, 2011).
For the design with $N$ distinct points for each $x_j$, the
similar rate
$\max\left(\sum_{j \in \Jc} N^{-r_j}, \frac{d_0 \ln(d/d_0)}{N}\right)$,
where $r_j=2s_j/(2s_j+1)$,
was derived in Raskutti, Wainright \& Yu (2012).
\ignore{
We argue that when it does not hold, trivial zero estimators
$\hat{f}_j=0,\;j=1,\ldots,d$ achieve a faster rate. The case $\ln(d/d_0)>N$
can be viewed as a {\em super}-sparse setup, where the proportion of nonzero
$f_j$ decreases exponentially with $N$.
}

We also propose a rate-optimal estimator for estimating
sparse additive models (\ref{eq:spad}) which is adaptive to the unknown
parameters $(s_j,R_j),\;j \in \Jc$ of Sobolev ellipsoids and to the unknown
sparsity $d_0$. The estimation is performed in the Fourier domain and
is based on
identifying nonzero vectors of (univariate) discrete Fourier coefficients
$c_j$ by imposing a penalty on the number of nonzero $c_j$'s and estimating their
components by truncating the corresponding series of empirical Fourier
coefficients of the data,
and can be efficiently computed.
The resulting estimator is developed within a
Bayesian framework and can be viewed as a maximum {\em a posteriori} (MAP)
sparse additive estimator. From a frequentist view, it corresponds to
penalized maximum likelihood estimation of $c_j$ with the complexity type of penalties
on the number of nonzero $c_j$ and numbers of their nonzero entries.

We compare the sparse additive MAP estimator with several existing
counterparts proposed recently in the literature, e.g., COSSO of Lin \& Zhang
(2006), SPAM of Ravikumar {\em et al.} (2009), sparse additive estimator
of Meier, van de Geer \&
B\"uhlmann (2009) and $M$-estimator of Raskutti, Wainwright \& Yu (2012)
(see also Koltchinskii \& Yuan, 2010 and Suzuki \& Sugiyama, 2013). In the
Fourier domain, the above estimators also correspond to penalized
maximum likelihood estimation of $c_j$ but with penalties on the
magnitudes of $c_{kj}$ rather than on their cardinality. However,
only the $M$-estimator is proved to
be rate-optimal (in the minimax sense) for the case when there are
$N$ distinct observations for each predictor $x_j$. Moreover, all
those procedures (except SPAM) are
not adaptive to the smoothness $s_j$ of $f_j$.

The paper is organized as follows. In Section \ref{sec:map} we derive the
sparse additive MAP estimator. Its asymptotic adaptive minimaxity is
established in Section \ref{sec:minimax}, where we compare it also with
its existing counterparts. The results of a simulation study are given
in Section \ref{sec:simul}. Some concluding remarks and possible extensions
are discussed in Section \ref{sec:disc}. All the proofs are placed in the
Appendix.

\section{MAP estimator} \label{sec:map}
\subsection{Main idea} \label{subsec:idea}

For any fixed $j=1,\ldots,d$, averaging a general additive model
(\ref{eq:admodel}) over all $N/n_j$ observations at points $i_j/n_j$
and using the identifiability conditions yields
\be \label{eq:barj}
\begin{split}
\bar{y}_j(i_j/n_j)=&\frac{n_j}{N}\sum_{i_1=0}^{n_1-1}\ldots\sum_{i_{j-1}=0}^{n_{j-1}-1}
\sum_{i_{j+1}=0}^{n_{j+1}-1}\ldots\sum_{i_d=0}^{n_d-1}y(i_1/n_1,\ldots,i_j/n_j,\ldots,i_d/n_d)\\
=&a_0+f_j(i_j/n_j)+\epsilon'(i_j/n_j), \quad i_j=0,\ldots,n_j
\end{split}
\ee
where $\epsilon'(i/n_j) \sim {\cal N}(0,\frac{n_j}{N}\sigma^2)$ and are
independent.

Equivalently, in the Fourier domain one has
\be \label{eq:fourier}
\xi_{kj}=c_{kj}+\frac{\sigma^2}{N} z_{kj},\quad
k=-(n_j-1)/2,\ldots,(n_j-1)/2;\;j=1,\ldots,d, \ee
where
$$
\xi_{kj}=\frac{1}{n_j}\sum_{i=0}^{n_j-1}\bar{y}_j(i/n_j)e^{\frac{2\pi I k i}{n_j}}
$$
are discrete (one-dimensional) Fourier coefficients of the vector $\bar{y}_j$,
$c_{kj}$ are given in (\ref{eq:disfourier}) and
$z_{kj}$ are independent standard complex normal variates.

The goal now is to
estimate the unknown discrete Fourier coefficients $c_{kj}$
in (\ref{eq:fourier}) by some $\hat{c}_{kj}$.
The resulting estimator
$\hat{f}$ in the original domain will then be
$$
\hat{f}(i_1/n_1,\ldots,i_d/n_d)=\hat{a}_0+\sum_{j=1}^d
\hat{f}_j(i_j/n_j)=\hat{a}_0+\sum_{j=1}^d \sum_{k=-(n_j-1)/2}^{(n_j-1)/2}
\hat{c}_{kj}e^{\frac{I 2 \pi k i_j}{n_j}}.
$$
Additivity of $f$ and Parseval's equality imply
$$
AMSE(\hat{f},f)=E|\hat{a}_0-a_0|^2+\sum_{j=1}^d
E||\hat{c}_j-c_j||^2_{n_j}
$$
and the original dimensionality of the problem $N$ is thus reduced to
$\sum_{j=1}^d (n_j-1)+1$ in the Fourier domain (recall that $c_{0j}=0$ for all $j$).

Estimate the overall mean $a_0$ by the overall sample mean $\bar{y}$.
Due the identifiability conditions $\sum_{i=0}^{n-1}f_j(i/n)=0$, we have
$$
\bar{y}=a_0+\epsilon^*,
$$
where $\epsilon^* \sim {\cal N}(0,\frac{\sigma^2}{N})$,
yielding $E|\bar{y}-a_0|^2=\frac{\sigma^2}{N}$.
Furthermore, we naturally set $\hat{c}_{0j}=0$ for all $j$ with no error and,
therefore, $\sum_{i=0}^{n_j-1}\hat{f}_j(i/n_j)=0$.

Recall now that we consider a {\em sparse} additive model
(\ref{eq:spad}), where most $f_j$ and, therefore, $c_j$ are zeros.
Under the assumption $f_j \in \Fj,\;j \in \Jc$, the corresponding $c_{kj}$
decrease polynomially in $k$ and $c_j$ can be well-approximated by several
first $c_{kj}$.
The proposed algorithm tries first to identify the set $\Jc$ of
nonzero vectors $c_j$ and then estimates their entries by truncating the
corresponding vectors $\xi_j$ of empirical discrete Fourier coefficients
in (\ref{eq:fourier}) at the properly adaptively chosen cut-points.

\subsection{Derivation} \label{subsec:deriv}
For nonzero vectors $c_j$ in (\ref{eq:fourier}) we consider truncated estimators
of the form $\hat{c}_{kj}=\xi_{kj},\;
|k|=1,\ldots,k_j$ and zero otherwise. Thus, if we knew the set of indices $\Jc$
of nonzero
$c_j$ and the cut-points $k_j,\;j\in\Jc$, we would
estimate $c_{kj},\;|k|=1,\ldots,k_j,\;j \in \Jc$ by the corresponding
$\xi_{kj}$ and set the others to zero. Since in reality they are unknown we
should estimate them from the data.

We use a Bayesian framework. Consider the following hierarchical
prior model on vectors $c_j$.
Let $d_{0}=|\Jc|=\#\{\, j\,:\, c_j \neq 0,\, j=1,\ldots,d\}$
be the number of nonzero $c_j$, and assume some prior
distribution $\pi(d_0)>0,\;d_0=0,\ldots,d$ on $d_0$.
For a given $d_{0}$, assume that all possible sets $\Jc$ of nonzero $c_j$
with $|\Jc|=d_0$
are equally likely, that is,
$$
P(\Jc \mid |\Jc|=d_{0})=\binom{d}{d_{0}}^{-1}.
$$
Obviously, $k_j|(j \in \Jo) \sim \delta(0)$ and, thus, $c_j|(j \in \Jo)
\sim \delta(0)$.
For nonzero $c_j$ we assume some
independent priors $\pi_j(k_j)|(j \in \Jc)>0,\;k_j=1,\ldots,(n_j-1)/2$.
To complete the prior we place independent normal priors for nonzero
$c_{kj}\sim {\cal N}(0,\gamma \frac{\sigma^2}{N}),\;j \in \Jc,\;|k|=1,\ldots,
k_j$, where $\gamma>0$. One can also consider different $\gamma_j$.

By a straightforward Bayesian calculus, the posterior probability of a given
set $\Jc$ and the corresponding $k_{j}$'s is
$$
P\left(\Jc;k_1,\ldots,k_{d_0}\,|\xi\right)\propto\pi_{0}(d_{0})\,
\binom{d}{d_{0}}^{-1}\prod_{j\in\Jc}\left\{\pi_{j}(k_j)
(1+\gamma)^{-k_j}\exp
{\left(\frac{\gamma}{1+\gamma}\frac{\sum_{|k|=1}^{k_j}
|\xi_{kj}|^2}{2\sigma^2/N}\right)}\right\}.
$$
Given the posterior distribution $P\left(\Jc;k_1,\ldots,k_{d_0}|\xi\right)$
we apply the maximum {\em a posteriori} (MAP) rule to find the most likely
set of nonzero vectors $\Jc$ and the corresponding cut-points $k_j,\;j\in \Jc$:
\begin{align}
\max_{\Jc;k_1,\ldots,k_{d_0}}&\left\{\sum_{j\in\Jc}\left[
\sum_{|k|=1}^{k_j}
|\xi_{kj}|^2+2\frac{\sigma^2}{N}\left(1+\frac{1}{\gamma}\right)\log\left(\pi_j(k_j)
(1+\gamma)^{-k_j}\right)\right] \right. \nonumber\\
&\left. +2\frac{\sigma^2}{N}\left(1+\frac{1}{\gamma}\right)
\log\left(\pi_0(d_0)\binom{d}{d_0}^{-1}\right)\right\}
\label{eq:map0}
\end{align}
To solve (\ref{eq:map0}), define $\hat{k}_{j}$ by
\begin{align}
\hat{k}_{j} & =\arg\min_{1\leq k_{j}\leq (n_j-1)/2}\left[\sum_{k:\,|k|>k_{j}}|\xi_{kj}|^2
+2\frac{\sigma^2}{N}\left(1+\frac{1}{\gamma}\right)\log\left(\pi_{j}^{-1}
\left(k_{j}\right)(1+\gamma)^{k_j}\right)\right]\nonumber \\
& =\arg\min_{1\leq k_j \leq (n_j-1)/2}\left[-\sum_{|k|=1}^{k_j}|\xi_{kj}|^{2}+
2\frac{\sigma^2}{N}\left(1+\frac{1}{\gamma}\right)\log\left(\pi_{j}^{-1}\left(k_{j}\right)(1+\gamma)^{k_j}\right)\right]\label{eq:map1}
\end{align}
for each $j=1,\ldots,d$.
The MAP rule in (\ref{eq:map0}) is then equivalent to minimizing
\begin{equation}
\sum_{j\in\Jc}\left\{-\sum_{|k|=1}^{\hat{k}_j}|\xi_{kj}|^2
+2\frac{\sigma^2}{N}\left(1+\frac{1}{\gamma}\right)\log\left(\pi_j^{-1}
(\hat{k}_j)(1+\gamma)^{\hat{k}_j}\right)+2\frac{\sigma^2}{N}
\left(1+\frac{1}{\gamma}\right)\log\left(\pi_0^{-1}\left(d_0\right)\binom{d}{d_0}\right)\right\}\label{eq:criterion}
\end{equation}
over all subsets of indices $\Jc\subseteq\{1,...,d\}$, where $d_0=|\Jc|$, and
the resulting algorithm for solving (\ref{eq:map0}) is then as follows:

\medskip
\noindent
{\em Algorithm}
\begin{enumerate}

\item
For each $j=1$ to $d$, find $\hat{k}_{j}$ in (\ref{eq:map1}) and calculate
$$
W_{j}=-\sum_{|k|=1}^{\hat{k}_j}|\xi_{kj}|^2+2\frac{\sigma^2}{N}\left(1+\frac{1}
{\gamma}\right)\log\left(\pi_j^{-1}(\hat{k}_j)(1+\gamma)^{\hat{k}_j}\right).
$$

\item
Order $W_{j}$ in ascending order $W_{(1)}\leq\ldots\leq W_{(d)}$
and find $\hat{d}_{0}$:
$$
\hat{d}_{0}=\arg\,\min_{0\leq d_{0}\leq d}\sum_{j=1}^{d_0}\left\{W_{(j)}+
2\frac{\sigma^2}{N}\left(1+\frac{1}{\gamma}\right)\log\left(\pi^{-1}\left(d_0\right)\binom{d}{d_0}\right)\right\}.
$$

\item
Let $\hat{\Jc}$ be the set of indices corresponding
to the $\hat{d}_{0}$ smallest $W_j$. Set $\hat{c}_{j}=0$
for all $j\in\hat{\Jo}$ and $\hat{c}_{kj}=
\xi_{kj}~\mathbb{I}\{1 \leq |k|\leq\hat{k}_{j}\},\;k=0,\ldots,n_j; j \in \hat{\Jc}$ (recall that due to the identifiability conditions, $\hat{c}_{0j}=0$ for all $j$).

\end{enumerate}

One can easily verify that the resulting MAP estimators $\hat{c}_j$ can be
equivalently viewed as penalized likelihood estimators of $c_j$ in (\ref{eq:fourier}) of the form
\begin{equation} \label{eq:MAP}
\min_{\tilde{c}_j,\ldots,\tilde{c}_d}\left\{\sum_{j=1}^d\left(||\xi_{j}-\tilde{c}_j||^{2}_{n_j}+Pen_j(k_j)\right)+Pen_0(d_0)\right\}
\end{equation}
with the complexity penalty
\begin{equation}
\label{eq:penalty0}
Pen_0(d_0)=2\frac{\sigma^2}{N}\left(1+\frac{1}{\gamma}\right)\log\left(\pi_{0}^{-1}(d_0)\binom{d}{d_0}\right)
\end{equation}
on the number of nonzero $\tilde{c}_j$ and the complexity penalties
\begin{equation}
\label{eq:penaltyj}
Pen_j(k_j)=2\frac{\sigma^2}{N}\left(1+\frac{1}{\gamma}\right)\log\left(\pi_j^{-1}\left(k_j\right)(1+\gamma)^{k_j}\right), \quad k_j=1,\ldots,(n_j-1)/2
\end{equation}
on the number of nonzero entries $2k_j$ of $\tilde{c}_j$.

\section{Theoretical properties} \label{sec:minimax}
\subsection{Upper bound} \label{subsec:upper}
In this section we establish theoretical properties of the proposed
sparse additive MAP estimator and establish its adaptive minimaxity with respect to the
$AMSE(\hat{f},f)=\sum_{j=1}^d AMSE(\hat{f}_j,f_j)$. As we have
mentioned, due to the Parseval's equality,
$AMSE(\hat{f},f)=\frac{\sigma^2}{N}+\sum_{j=1}^d
E||\hat{c}_j-c_j||^2_{n_j}$, where $\hat{c}_j$ and $c_j$ are
discrete Fourier coefficients of $\hat{f}_j$ and $f_j$ respectively
(see (\ref{eq:fourier})).

We start from a general upper bound on the $AMSE(\hat{f},f)$.
Recall that $N=\prod_{j=1}^d n_j$.
\begin{proposition}[general upper bound] \label{prop:upper}
Consider the sparse additive model (\ref{eq:spad}). Let $\hat{c}_1,\ldots,\hat{c}_d$
be the sparse additive MAP estimators
(\ref{eq:MAP}) of the Fourier coefficients vectors $c_1,\ldots,c_d$ in (\ref{eq:disfourier}) with the complexity penalties (\ref{eq:penalty0}) and (\ref{eq:penaltyj}).
Assume that
$\pi_j(k) \leq e^{-c(\gamma)k},\;k=1,...,(n_j-1)/2$ for all $j=1,\ldots,d$,
where $c(\gamma)=8(\gamma+3/4)^2>9/2$.
Then,
\be
\begin{split} \nonumber
AMSE(\hat{f},f) \leq C_1(\gamma) \min_{\Jo \subseteq \{1,...,d\}} & \left[
\sum_{j \in \Jc}
\min_{1 \leq k_j \leq (n_j-1)/2}
\left\{\sum_{|k| = k_j+1}^{(n_j-1)/2} |c_{kj}|^2+Pen_j(k_j)\right\} \right.  \\
& +  \left. \sum_{j \in \Jo}\sum_{k=-(n_j-1)/2}^{(n_j-1)/2}
|c_{kj}|^2+Pen_0(|\Jc|)\right] +C_2(\gamma) \frac{\sigma_2}{N} \{1-\pi_0(0)\},
\end{split}
\ee
where $C_1(\gamma)$ and $C_2(\gamma)$ depend only on $\gamma$.
\end{proposition}

Proposition \ref{prop:upper} holds without any regularity conditions on nonzero
$f_j$. Now we consider $f_j \in \Fj,\;j \in \Jc$:

\begin{theorem}[upper bound over $\Fj$] \label{th:upper}
Consider the model (\ref{eq:spad}), where $\Jc \neq \emptyset$. Assume that
$f_j \in \Fj$ for all $j \in \Jc$.

Let $\hat{c}_1,\ldots,\hat{c}_d$
be the sparse additive MAP estimators (\ref{eq:MAP}) of the Fourier coefficients
vectors $c_1,\ldots,c_d$ in (\ref{eq:disfourier}) with the complexity penalties
(\ref{eq:penaltyj})--(\ref{eq:penalty0}). Assume that
there exist constants $C_0, C_1>0$ such that
\begin{enumerate}
\item $\pi_0(h) \geq (h/d)^{C_0 h},\; h=1,\ldots,\lfloor d/e\rfloor$ and
$\pi_0(d) \geq e^{-C_0 d}\;;$
\item $e^{-C_1 k} \leq \pi_j(k) \leq e^{-c(\gamma)k},\;k=1,\ldots,(n_j-1)/2,\;j=1,\ldots,d$
\end{enumerate}

Then, for any $\Jc \subseteq \{1,\ldots,d\}$ with $|\Jc|=d_0$ and all
$\Fj,\; j \in \Jc$,
\be \label{eq:upper}
\sup_{f_j \in \Fj,
j \in \Jc}AMSE(\hat{f},f)
\leq C_1(\gamma) \max\left\{\sum_{j \in \Jc}
\min\left(N^{-\frac{2s_j}{2s_j+1}},\frac{n_j}{N}\right), \frac{d_0 \ln(d/d_0)}{N}\right\},
\ee
where $C_1(\gamma)$ is some constant  depending only on $\gamma$.
\end{theorem}

One can easily verify that the conditions on priors $\pi(\cdot)$ and
$\pi_j(\cdot)$ required in Theorem \ref{th:upper} are satisfied for
the (truncated) geometric priors $\pi_0(h) \propto
q^h,\;h=1,\ldots,d$ and $\pi_j(k) \propto q_j^k,\;k=1,\ldots,(n_j-1)/2$ for
some $0 < q, q_j < 1$ corresponding respectively to the complexity
penalties $Pen_0(h)\sim2C(\gamma)\frac{\sigma^2}{N}h(\ln(d/h)+1)$ of
the $2h\ln(d/h)$-type and the AIC type $Pen_j(k) \sim 2
C(\gamma)\frac{\sigma^2}{N}k$ for some $C(\gamma)>1$.

\subsection{Asymptotic minimaxity} \label{subsec:lower}
To assess the goodness of the upper bound for the AMSE of the MAP
estimator established in Theorem \ref{th:upper} we derive the corresponding
minimax lower bounds.

We start from the following proposition establishing the
minimax lower bound for estimating a single $f_j \in \Fj$ in the model (\ref{eq:barj}):
\begin{proposition}[minimax lower bound for a single $f_j \in \Fj$]
\label{prop:lower}
Consider the model (\ref{eq:barj}), where  $f_j \in \Fj$.
There exists a constant $C_2>0$ such that
$$
\inf_{\tilde{f}_j} \sup_{f_j \in \Fj} AMSE(\tilde{f}_j,f_j) \geq
C_2 \min\left(N^{-\frac{2s_j}{2s_j+1}},\frac{n_j}{N}\right),
$$
where the infimum is taken over all estimators $\tilde{f}_j$ of $f_j$.
\end{proposition}

We now use this result to obtain the minimax lower bound for the AMSE in estimating
$f$ in the sparse additive model (\ref{eq:spad}):
\begin{theorem}[minimax lower bound] \label{th:lower}
Consider the model (\ref{eq:spad}), where $f_j \in \Fj,\;j \in \Jc$.
There exists a constant $C_2>0$ such that
\be \label{eq:lower}
\inf_{\tilde{f}}\sup_{f_j \in \Fj,
j \in \Jc} AMSE(\hat{f},f) \geq
C_2 \max\left\{\sum_{j\in\Jc}\min\left(N^{-\frac{2s_j}{2s_j+1}},\frac{n_j}{N}\right),\frac{d_0\ln(d/d_0)}{N}\right\},
\ee
where the infimum is taken over all estimators $\tilde{f}$ of $f$.
\end{theorem}

Theorems \ref{th:upper} and \ref{th:lower} shows that as both the
sample sizes $n_j$'s and the dimensionality $d$ increase, the
asymptotic minimax convergence rate is either of order
$\sum_{j\in\Jc}\min\left(N^{-\frac{2s_j}{2s_j+1}},\frac{n_j}{N}\right)$
or $N^{-1}d_0\ln(d/d_0)$. The former corresponds to the optimal
rates of estimating $d_0$ single $f_j \in \Fj$, while the latter is
due to error in selecting a subset of $d_0$ nonzero $f_j$ out of $d$
and commonly appears in various related model selection setups (see,
e.g., Abramovich \& Grinshtein, 2010, 2013; Raskutti, Wainwright \&
Yu, 2011, 2012; Rigollet \& Tsybakov, 2011). Dominating term depends
on the smoothness of $f_j$'s (relatively to the sample sizes
$n_j$'s) and sparsity of the problem.

Furthermore, the proposed sparse additive MAP estimator with the
priors $\pi_0(\cdot)$ and $\pi_j(\cdot)$ corresponding to
$2d_0\ln(d/d_0)$-type and AIC-type penalties respectively is
simultaneously minimax rate-optimal over the entire range of sparse
and dense amalgams of Sobolev balls $\Fj$. \ignore{ excluding
super-sparse cases, where $\ln(d/d_0)>NR_j^2$ for $j \in \Jc$. When
such super-sparsity presents, the minimax lower bound
(\ref{eq:lower}) can be reduced by the smoothness of $f_j$ trivial
zero estimators $\tilde{f}_j=0,\;j=1,\ldots,d$. Indeed, in this case
evidently
$$
\sup_{f_j: f_j \in \Fj} AMSE(\tilde{f},f)=
\sum_{j \in \Jc}\sup_{c_j \in \Tj} ||c_j||^2_n.
$$
The least favorable vector of discrete Fourier coefficients $c_j \in \Tj$
with the maximal $l^2_n$-norm $||c_j||^2_n$
is the one with only two nonzero entries $c_{-1j}=c_{1j}=R_j/\sqrt{2}$.
Hence,
$\sum_{j \in \Jc}\sup_{c_j \in \Tj}||c_j||^2_n=\sum_{j \in \Jc} R_j^2<
d_0\ln(d/d_0)/N$
and the resulting $AMSE(\tilde{f},f)$ is of order less than (\ref{eq:lower}).
}

\subsection{Comparison with other existing estimators} \label{subsec:comp}
As we have already mentioned, various estimators for the sparse additive
model (\ref{eq:spad}) have been recently proposed in the literature.
It can be shown that being adapted to the considered setup,
they can be also equivalently formulated in the Fourier
domain as penalized maximum likelihood estimators of $c_j$ but with
penalties on the magnitudes of $c_{kj}$ rather than
complexity-type penalties as for the proposed sparse additive MAP estimator.

Thus, the additive COSSO method of
Lin \& Zhang (2006, Section 4) in this case can be written as
\be \label{eq:cosso}
\arg \min_{\tilde{c}_j,\ldots,\tilde{c}_d;~\theta_1>0,\ldots,\theta_d>0}
\left\{\sum_{j=1}^d||\xi_{j}-\tilde{c}_j||^{2}_{n_j}+\sum_{j=1}^d \theta_j^{-1} \sum_{k=-(n_j-1)/2}^{(n_j-1)/2}|k|^{2s_j}
|\tilde{c}_{kj}|^2+\lambda \sum_{j=1}^d \theta_j \right\}.
\ee
The form of the estimator (\ref{eq:cosso}) is very
similar to the common spline smoothing which is equivalent to linear shrinkage
in the Fourier domain (e.g., Wahba, 1990) with smoothing
parameters $\theta_j$ but with the additional penalty on their sum. The latter
makes the set of optimal $\theta_j$ to be sparse and, therefore, yields
zero components $\hat{c}_j$ in the resulting COSSO estimators. To the best
of our knowledge, there are no results on the convergence rates for the COSSO.

Similarly, the sparse additive estimator of Meier, van de Geer \& B\"uhlmann (2009)
can be presented as
\be \label{eq:meier}
\arg \min_{\tilde{c}_j,\ldots,\tilde{c}_d}
\left\{\sum_{j=1}^d||\xi_{j}-\tilde{c}_j||^{2}_{n_j}
+\lambda_1 \sum_{j=1}^d\sqrt{||\tilde{c}_j||^2_{n_j}+\lambda_2
\sum_{k=-(n_j-1)/2}^{(n_j-1)/2}
|k|^{2s_j}|\tilde{c}_{kj}|^2} \right\} ,
\ee
where penalizing $||\tilde{c}_j||_{n_j}$ encourages sparsity, while
the additional penalty term controls the smoothness of the estimators.
For $N$ distinct observations for each $x_j$, from the results of
Meier, van de Geer \& B\"uhlmann (2009, Remark 2) it
follows that their estimator has a sub-optimal rate
$O\left(\sum_{j \in \Jc} \left(\frac{\ln d}{N}\right)^{\frac{2s_j}{2s_j+1}}\right)$.

Applied to $f_j \in \Fj$, a regularized $M$-estimator of Raskutti,
Wainwright \& Yu (2012) is \be \label{eq:m}
\arg \min_{\tilde{c}_j,\ldots,\tilde{c}_d}
\left\{\sum_{j=1}^d||\xi_{j}-\tilde{c}_j||^{2}_{n_j} +\lambda_1
\sum_{j=1}^d ||\tilde{c}_j||_{n_j} + \lambda_2 \sum_{j=1}^d
\sqrt{\sum_{k=-(n_j-1)/2}^{(n_j-1)/2} |k|^{2s_j}|\tilde{c}_{kj}|^2} \right\}
\ee which is similar to (\ref{eq:meier}) but separates the penalties
on sparsity and smoothness into two additive terms. For the design
with $N$ distinct observations for each $x_j$, the estimator
(\ref{eq:m}) achieves the minimax rate $O\left(\min\left(\sum_{j \in
\Jc} N^{-\frac{2s_j}{2s_j+1}}\right), \frac{d_0
\ln(d/d_0)}{N}\right)$. Similar results for the $M$-estimator (\ref{eq:m})
were obtained in Koltchinskii \& Yuan (2011) and Suzuki \& Sugiyama (2013)
under some additional conditions.

The serious disadvantage of all the above estimators is that they
are defined for penalties involving $s_j$ and, hence, are inherently
not adaptive to the smoothness of $f_j$ which can rarely be assumed known.

The SPAM estimator of Ravikumar {\em et al.} (2009) for
the considered setup becomes
\be \label{eq:spam}
\arg \min_{\tilde{c}_j,\ldots,\tilde{c}_d}
\left\{\sum_{j=1}^d||\xi_{j}-\tilde{c}_j||^2_{n_j}+
\lambda \sum_{j=1}^d \sqrt{2k_j}||\tilde{c}_j||_{k_j} \right\}
\ee
for the fixed truncation cut-points $k_j$. In this form, SPAM is closely related
to the group lasso estimator of Yuan \& Lin (2006) and can be obtained explicitly:
\be \label{eq:grouplasso}
\hat{c}_j=\left(1-\frac{(\lambda/2) \sqrt{2k_j}}{||\tilde{\xi}_j||_{k_j}}\right)_+ \tilde{\xi}_j,
\ee
where $\tilde{\xi}_j$ is $\xi_j$ truncated at $k_j$. Ravikumar {\em et al.} (2009)
show persistency of their estimator but do not provide results on
convergence rates of its AMSE.

Finally, we can mention Guedj \& Alquier (2013) that
considered a Bayesian model similar to that proposed in this paper with
geometric priors $\pi_0(\cdot)$ and $\pi_j(\cdot)$. They estimated $c_j$ by the
corresponding posterior means and for the case of $N$ distinct
observations for each $x_j$, showed that the resulting estimator is
asymptotically
nearly-minimax (up to an additional log-factor) over Sobolev classes.
A similar Bayesian estimator of Suzuki (2012) achieves the optimal rate but
for smaller functional classes.
The practical implementation of these procedures involves however
high-dimensional MCMC algorithms.

\section{Simulation study} \label{sec:simul}
To illustrate the performance of the proposed sparse additive MAP estimator
we conducted a simulation study. Similar to Example 1 of Lin \& Zhang (2006), Example 3 of Meier, van de Geer \& B\"uhlmann (2009) and Example 3 of
Guedj \& Alquier (2013), we considered the sparse additive model
(\ref{eq:spad}) with $d=50$ and four nonzero components $f_j$ ($d_0=4$):
\begin{eqnarray*}
& f_1(x)=x \\
& f_2(x)=(2x-1)^{2} \\
& f_3(x)=\frac{\sin(2\pi x)}{2-\sin(2\pi x)}\\
& f_4(x)=0.1\sin(2\pi x)+0.2\cos(2\pi x)+0.3\sin^{2}(2\pi x)+0.4\cos^{3}(2\pi x)+0.5\sin^{3}(2\pi x)
\end{eqnarray*}
but on the regular lattice $[0,1]^{50}$.
We used $n=101$ and, therefore, $N=101^{50}$.
Each nonzero $f_j$ was standardized to have
\begin{eqnarray*}
& \frac{1}{n}\sum_{i=0}^{n-1}f_{j}(i/n)=0, \\
&\frac{1}{n}\sum_{i=0}^{n-1}f_{j}^{2}(i/n)=1.
\end{eqnarray*}

The noisy data was generated according to model (\ref{eq:barj}) by
adding independent random Gaussian variates ${\cal
N}(0,\frac{n}{N}\sigma^2)$ to
$f_j(i/n),\;i=0,\ldots,n-1;~j=1,\ldots, d$. The values of the noise
variance $\sigma^2$ were chosen to correspond to values 1, 5 and 10
for the signal-to-noise ratio (SNR) defined as
$SNR=Var(f_j)/(\frac{n}{N}\sigma^2)= \frac{N}{\sigma^2 n}$.
Performing the discrete Fourier transform of the noisy data yielded
the equivalent model (\ref{eq:fourier}) in the Fourier domain. We
applied then the proposed MAP algorithm to corresponding noisy
Fourier coefficients $\xi_{kj}$ using truncated geometric priors for
$\pi_0(\cdot)$ and $\pi_j(\cdot)$ with $q=q_j=0.5$ and $\gamma=5$.
The noise level $\sigma$ was assumed unknown and estimated from the
data. Since the vector of the true Fourier coefficients $c_j$ in 
(\ref{eq:fourier}) lies in a Sobolev ellipsoid, the sequence
$|c_{kj}|$ decays to zero polynomially with $k$. Thus, for large
$k$, the empirical Fourier coefficients $\xi_{kj}$ in
(\ref{eq:fourier}) are mostly pure noise. To correct for the bias
due to the possible presence of several large coefficients, we
robustly estimated $\sigma/\sqrt{N}$ from $\xi_{kj}$ for large $k$ as follows:
$$
\frac{\hat{\sigma}}{\sqrt{N}}=\frac{{\sqrt 2}~MAD\left(\{Re(\xi_{kj}),Im(\xi_{kj})\},\;
\;k=0.8~\frac{n_j-1}{2},\ldots, \frac{n_j-1}{2};\;j=1,\ldots,50\right)}{0.6745}~.
$$
This is similar to a standard practice for estimating $\sigma$ from wavelet
coefficients at
the finest resolution level in wavelet-based methods (see, e.g., Donoho
\& Johnstone, 1994). The resulting estimates
for $\sigma$ were very precise for all SNRs.

We compared also the resulting sparse additive MAP estimator with
the SPAM estimator (\ref{eq:spam}) of Ravikumar {\em et al.} (2009)
which for the considered model is essentially the group lasso
estimator of Yuan \& Lin (2006) and is available in the closed form
in the Fourier domain -- see (\ref{eq:grouplasso}).
For the SPAM estimator we used the
same cut-points $\hat{k}_j$ from (\ref{eq:map1}) as for the MAP, and
the oracle chosen threshold $\lambda$ that minimizes the
$AMSE(f,\hat{f}^{SPAM})= \sum_{j=1}^d
||\hat{c}^{SPAM}_j(\lambda)-c_j||^2_n$ estimated by averaging over a
series of 1000 replications for each value of $\lambda$ using a grid
search. The resulting choices were $\lambda=0.26$ for $SNR=1$,
$\lambda=0.10$ for $SNR=5$ and $\lambda=0.06$ for $SNR=10$. Thus,
the oracle $\lambda$ decreased with increasing SNR.

\begin{figure}[h]
\begin{center}
\includegraphics[height=12.5cm]{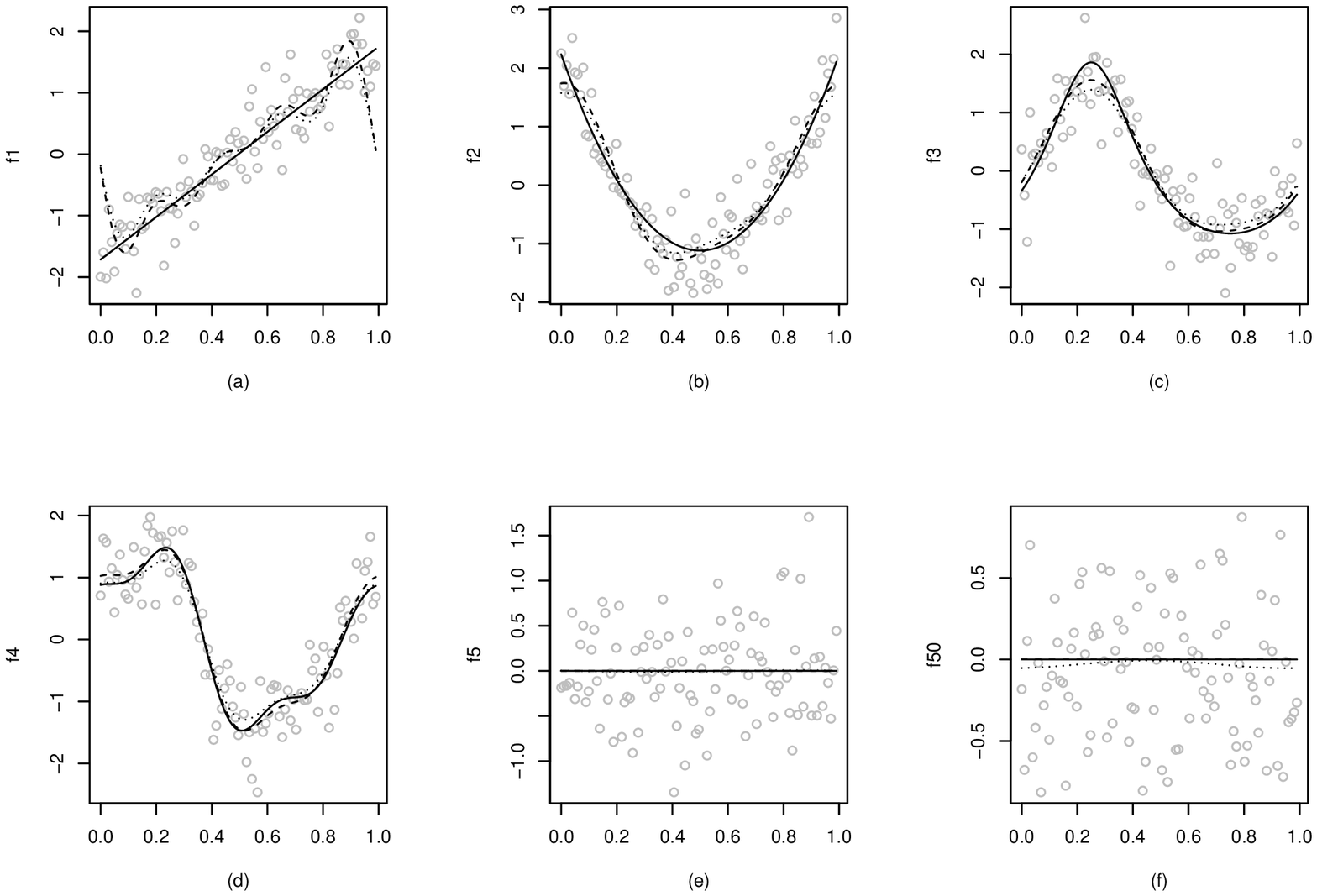} 
\caption{Examples of MAP (dashed lines) and SPAM (dotted lines) estimators
for various $f_j$ (solid lines): $f_1$ (a), (b) $f_2$ (b), $f_3$ (c),  $f_4$ (d)
and two zero $f_j$ (e)-(f) ($SNR=5$).}
\label{fig:example}
\end{center}
\end{figure}

For each SNR level we calculated the (global) $AMSE$ for both methods and
analyzed also their performance for each individual $f_j$. Thus, $AMSE_1,\;AMSE_2\;
AMSE_3\;AMSE_4$ are the AMSEs for the corresponding four nonzero
$f_j,\;j=1,\ldots 4$, while $AMSE_0$ is the average AMSE over all 46 zero $f_j$.
In addition, we compared the two methods for identifying nonzero
$f_j$ though it is a somewhat different problem from our original goal of
estimating functions in quadratic norm and calculated $\hat{d}_0=\#\{j:
\hat{f}_j \ne 0,\;j=1,\ldots,50\}$. The results are summarized in
Table \ref{tab:amse} below. See also Figure \ref{fig:boxplot} for the
corresponding boxplots. Figure \ref{fig:example} gives typical examples of
estimators obtained by both methods for nonzero and zero $f_j$.

\begin{table}[htb]
\begin{center}
\caption{AMSE averaged over 1000 replications for various SNR.
}{%
\begin{tabular}{c|l|c|ccccc|c}
\hline
SNR  &  method & $AMSE$  & $AMSE_1$ & $AMSE_2$ & $AMSE_3$ & $AMSE_4$ & $AMSE_0$ & $\hat{d}_0$ \\
\hline

1    & MAP             & 0.6242  & 0.3083   & 0.1023 & 0.0926 & 0.1209 & 0.0000 & 4.0\\
     &SPAM($\lambda=0.26$) & 0.8007 & 0.3371 & 0.1283 & 0.1178 & 0.1467 & 0.0015 & 19.3       \\
\hline
5    & MAP   & 0.1937 & 0.1334 & 0.0285 & 0.0157 & 0.0161 & 0.0000 & 4.0     \\
     &SPAM($\lambda=0.10$) & 0.2632 & 0.1492 & 0.0373 & 0.0238 & 0.0282 & 0.0005 & 25.7    \\
\hline
10   & MAP             &  0.1285 & 0.0936 & 0.0182 & 0.0099 & 0.0067 & 0.0000 & 4.0      \\
     &SPAM($\lambda=0.06$) & 0.1686 & 0.1021 & 0.0220 & 0.0131 & 0.0114 & 0.0004& 32.3     \\
\hline

\end{tabular}}
\label{tab:amse}
\end{center}
\end{table}

\begin{figure}[h]
\begin{center}
\includegraphics[height=6.2cm]{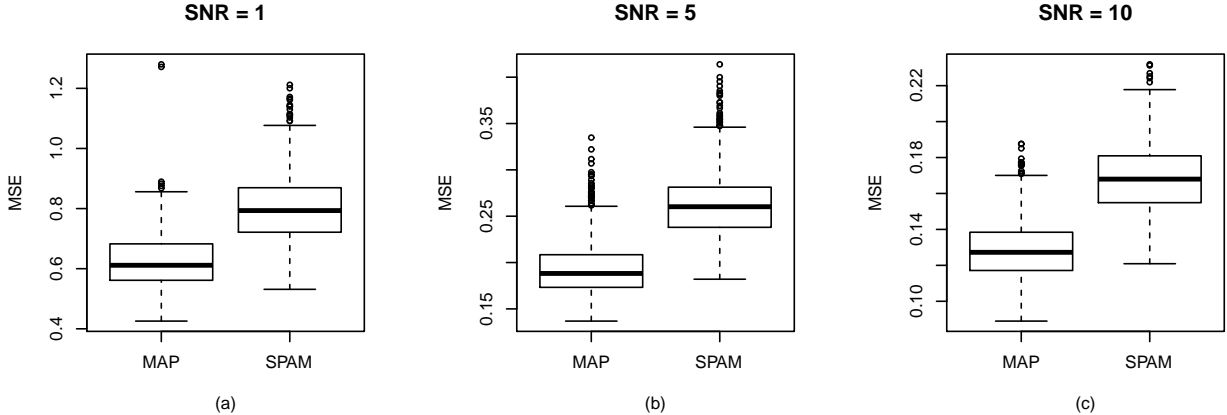} 
\caption{Boxplots for (global) AMSE for various SNR.}
\label{fig:boxplot}
\end{center}
\end{figure}

The results in Table \ref{tab:amse} show that MAP consistently outperforms
SPAM (even with the oracle choices for $\lambda$)
both globally and for each individual component $f_j$. For both
methods the main contribution to the global AMSE came from estimating nonzero
$f_j$. The MAP estimator almost perfectly identified the set of nonzero $f_j$
while the oracle
choices for $\lambda$ in SPAM were quite small and, as a result, too many
$\hat{f}_j$ were nonzero (see, e.g., Figure \ref{fig:example} (f)). In fact, it is a known common phenomenon for lasso-type
estimators.

\section{Concluding remarks} \label{sec:disc}
We considered sparse additive regression on a regular lattice, where
the univariate components $f_j$ of the unknown response function $f$
belong to Sobolev balls. We established the minimax convergence
rates of estimating $f$ and proposed an adaptive Fourier-based
estimator which is rate-optimal over the entire range of Sobolev
classes of different sparsity and smoothness. The resulting
estimator was developed within Bayesian formalism but can also be
viewed, in fact, as a penalized maximum likelihood estimator of the
Fourier coefficients of $f$ with certain complexity penalties on the
number of nonzero $f_j$ and on the numbers of nonzero entries of
their Fourier coefficients $c_j$. It can be efficiently computed and
the presented simulation study demonstrates its good
performance.

The results of the paper can be extended to more general Besov classes of
functions using the wavelet series expansions of $f_j$.
The corresponding vectors of wavelet coefficients will lie then within weak
$l_p$-balls (e.g., Johnstone, 2013, Section 9.7) and one can apply the results of
Abramovich \& Grinshtein (2013) for estimating a sparse group of sparse
vectors from weak $l_p$-balls. The extension is quite straightforward though
the details should be worked out. In particular, the resulting MAP estimator
should mimic (hard) thresholding within each nonzero vector of wavelet
coefficients
instead of truncation as in the considered case of Fourier series
(see Abramovich \& Grinshtein, 2013).

\section*{Acknowledgement} The work was supported by the Israel Science Foundation
(ISF), grant ISF-820/13. We are grateful to Anestis Antoniadis, Alexander
Goldenshluger and Vadim Grinshtein for fruitful discussions and valuable remarks.
Helpful comments by the Editor and an anonymous referee are gratefully acknowledged.

\section*{Appendix}
Throughout the proofs we use $C$ to denote a generic positive constant, not
necessarily the same each time it is used, even within a single equation.
Similarly, $C(\gamma)$ is a generic positive constant depending on $\gamma$.

\subsection*{Proof of Proposition \ref{prop:upper}} \label{subsec:proofpropupper}

As we have mentioned before,
the proposed sparse additive MAP estimator (\ref{eq:MAP})
can be equivalently viewed a penalized maximum
likelihood estimator with complexity penalties (\ref{eq:penalty0}) and (\ref{eq:penaltyj}). We can apply then
the general results of Birge \& Massart (2001) for complexity
penalized estimators.

Rewrite first the model (\ref{eq:fourier}) in a different form.
Set $\xi=\left(\xi_{-(n_1-1)/2,1},\ldots,\xi_{(n_1-1)/2,1},\ldots,\xi_{-(n_d-1)/2,d},\ldots,\xi_{(n_d-1)/2,d}\right)^{t}$
to be an amalgamated vector of length
$N_0=\sum_{j=1}^d n_j$ of $d$ vectors  $\xi_1,\ldots,\xi_d$. Similarly, define $N_0$-dimensional amalgamated vectors
$c=\left(c_{-(n_1-1)/2,1},\ldots,c_{(n_1-1)/2,1},\ldots,c_{-(n_d-1)/2,d},\ldots,c_{(n_d-1)/2,d}\right)^{t}$ and
$z=\left(z_{-(n_1-1)/2,1},\ldots,z_{(n_1-1)/2,1},\ldots,z_{-(n_d-1)/2,d},\ldots,z_{(n_d-1)/2,d}\right)^{t}$.
The original model (\ref{eq:fourier}) can be rewritten then as
\begin{equation}
\xi_{i}=c_{i}+\frac{\sigma^2}{N}z_{i}, \quad i=1,\ldots,N_0 \label{eq:model}
\end{equation}
where $z_i$ are independent standard complex normal variates.
Define an indicator vector
$v$ by $v_{i}=\mathbb{\mathbb{I}}\left\{ c_{i}\neq0\right\}$, $i=1,\ldots,N_0$.
Thus, in terms of model (\ref{eq:model}), $k_{j}=(1/2)\sum_{i=S_{j-1}+1}^{S_j}v_{i}$,
where $S_j=\sum_{l=1}^{j-1}n_l$, and $d_0=\#\{j: k_j>0\}$.
For a given $v$, let $D_{v}=2\sum_{j=1}^{d}k_{j}=\#\left\{ i:\, v_{i}=1,\;
i=1,\ldots,N_0\right\}$ be the overall number of nonzero entries of $c$, and
define
\begin{equation}
L_{v}=\begin{cases}
\frac{1}{D_{v}}\left\{\sum_{j=1}^{d}\log\left(\pi_{j}^{-1}\left(k_{j}\right)\right)+\log\left(\pi_{0}^{-1}\left(d_{0}\right)\binom{d}{d_{0}}\right)\right\} & \mbox{if}\: v\neq0\\
\\
\log\pi_{0}^{-1}(0) & \mbox{if}\: v=0.
\end{cases}\label{eq:weights}
\end{equation}
In the above notations the sparse additive MAP estimator
$\hat{c}=\left(\hat{c}_{-(n_1-1)/2,1},\ldots,\hat{c}_{(n_1-1)/2,1},\ldots,\hat{c}_{-(n_d-1)/2,d},\ldots,\hat{c}_{(n_d-1)/2,d}\right)^{t}$
is the penalized maximum likelihood estimator of $c$ with the complexity
penalty
\begin{eqnarray}
Pen(v) & = & 2\frac{\sigma^{2}}{N}\left(1+\frac{1}{\gamma}\right)
\left\{\sum_{j=1}^{d}\log\left(\pi_{j}^{-1}\left(k_{j}\right)\left(1+\gamma\right)^{k_j}\right)+\log\left(\pi_{0}^{-1}\left(d_{0}\right)\binom{d}{d_{0}}\right)\right\} \nonumber\\
& = & 2 \frac{\sigma^{2}}{N}\left(1+\frac{1}{\gamma}\right)D_{v}\left(L_{v}+\frac{1}{2}\log\left(1+\gamma\right)\right) \label{eq:penalty}
\end{eqnarray}
for $v \neq 0$, and $Pen(0)=2\frac{\sigma^{2}}{N}\left(1+\frac{1}{\gamma}\right)L_{0}$.

One can easily verify that
$$
\sum_{v\neq0}\exp\left\{ -D_{v}L_{v}\right\} = \sum_{k=1}^{d}\pi_{0}(k)=1-\pi_{0}(0).
$$

\ignore{
\begin{eqnarray*}
\sum_{v\neq0}\exp\left\{ -D_{v}L_{v}\right\}  & = & \sum_{v\neq0}\exp\left\{ -\sum_{j=1}^{d}\log\left(\pi_{j}^{-1}\left(k_{j}\right)\right)-\log\left(\pi_{0}^{-1}\left(d_{0}\right)\binom{d}{d_{0}}\right)\right\} \\
\\
 & = & \sum_{k=1}^{d}\left[\binom{d}{k}\exp\left\{ \log\left(\pi_{0}(k)\right)+\log\binom{d}{k}^{-1}\right\} \sum_{||v||_{0}=k}\exp\left\{ \sum_{j=1}^{d}\log\left(\pi_{j}\left(k_{j}\right)\right)\right\} \right]\\
\\
 & = & \sum_{k=1}^{d}\pi_{0}(k)=1-\pi_{0}(0).
\end{eqnarray*}
}
Furthermore, straightforward calculus similar to that in the
proof of Theorem 1 of Abramovich {\em et al.} (2007) implies that
under the conditions on the priors $\pi_j(\cdot)$ of Proposition
\ref{prop:upper}, the complexity penalty $Pen(v)$ in
(\ref{eq:penalty}) satisfies
$$
Pen(v)\geq C(\gamma) \frac{\sigma^{2}}{N}D_{v}\left(1+\sqrt{2L_{v}}\right)^{2},
\quad
$$
for some $C(\gamma)>1$.
One can then apply Theorem 2 of Birge \& Massart (2001) to have
\begin{eqnarray*}
\sum_{j=1}^{d}E\left(\Vert\hat{c}_{j}-c_{j}\Vert_{2}^{2}\right) &
\leq & c_{1}(\gamma)\min_{\Jo\subseteq\{1,\ldots,d\}}\left\{\sum_{j\in\Jc}
\min_{1\leq k_{j}\leq (n_j-1)/2}\left( \sum_{k:|k|>k_{j}}|c_{kj}|^{2}+Pen_{j}\left(k_{j}\right)\right) \right.
\\
&  &\left. +\sum_{j\in\Jo}\sum_{|k|=1}^{(n_j-1)/2}|c_{kj}|^{2}+Pen_{0}\left(d_{0}\right)\right\}+c_{2}(\gamma)\frac{\sigma^{2}}{N}\left(1-\pi_{0}(0)\right).
\end{eqnarray*}
Parseval's equality $AMSE(\hat{f},f)=\sum_{j=1}^{d}E\left(\Vert\hat{c}_{j}-c_{j}\Vert_{n_j}^{2}\right)+\frac{\sigma^2}{N}$ completes the proof.

\qed

\subsection*{Proof of Theorem \ref{th:upper}} \label{subsec:proofthupper}
Let ${\cal J}_0^{c*}$ be the true (unknown) subset of nonzero $c_j$ and
$d_0^*=|{\cal J}_0^{c*}|$.
Consider separately two cases.

\noindent
{\em Case 1}: $d_0^* \leq \lfloor d/e \rfloor$. Applying the general upper bound
established in Proposition \ref{prop:upper} for $\Jo={\cal J}_0^*$ yields
\be
\begin{split}
AMSE(\hat{f},f) & \leq C_1(\gamma) \left\{
\sum_{j \in {\cal J}_0^{c*}}
\min_{1 \leq k_j \leq (n_j-1)/2}
\left\{\sum_{|k| = k_j+1}^{(n_j-1)/2} |c_{kj}|^2+Pen_j(k_j)\right\}
+Pen_0(d_0^*)\right\} \\
& +C_2(\gamma) \frac{\sigma_2}{N} \{1-\pi_0(0)\}. \label{eq:case1}
\end{split}
\ee
Choose the cut-points $k_j=\left\lfloor\frac{1}{2}\min(N^{\frac{1}{2s_j+1}},n_j-1)\right\rfloor$ for $j \in {\cal J}_0^{c*}$.
If $k_j < (n_j-1)/2$, for $c_j \in \Tj$ we have
$\sum_{|k| = k_j+1}^{(n_j-1)/2} |c_{kj}|^2=O(k_j^{-2s_j})=O\left(N^{-\frac{2s_j}{2s_j+1}}\right)$, while for $k=(n_j-1)/2$, this term obviously disappears.
Furthermore, under the conditions on the priors $\pi_j(\cdot)$, the
corresponding penalties $Pen_j(\cdot)$ in (\ref{eq:penaltyj}) are of the AIC-type, where
$Pen_j(k_j) \sim 2C(\gamma)\frac{\sigma^2}{N}k_j=
O\left(\min\left(N^{-\frac{2s_j}{2s_j+1}},\frac{n_j}{N}\right)\right)$.
Hence, the first term $\sum_{j \in {\cal J}_0^{c*}}$ in the RHS of
(\ref{eq:case1}) is of the order $\sum_{j \in {\cal J}_0^{c*}}
\min\left(N^{-\frac{2s_j}{2s_j+1}},\frac{n_j}{N}\right)$.

Finally, $\binom{d}{d^*_0} \leq
\left(\frac{d}{d^*_0}\right)^{2d^*_0}$ for $d^*_0 \leq \lfloor d/e \rfloor$
(see, e.g. Lemma A1 of Abramovich {\em et al.}, 2010) and, therefore, the
conditions on $\pi_0(\cdot)$ imply
$$
Pen_0(d^*_0) \leq C(\gamma)\frac{\sigma^2}{N} d_0^* \log(d/d_0^*).
$$

\noindent
{\em Case 2: $\lfloor d/e \rfloor < d^*_0 \leq d$}.
In this case we apply Proposition \ref{prop:upper} for $\Jo=\emptyset$.
Evidently, $|\Jc|=d$ and $\Jc={\cal J}_0^* \bigcup {\cal J}_0^{c*}$.
Choose the cut-points
$k_j=\left\lfloor\frac{1}{2}\min(N^{\frac{1}{2s_j+1}},n_j-1)\right\rfloor$ for $j \in {\cal J}_0^{c*}$ as before and $k_j=1$ for $j \in {\cal J}_0^*$.
Then,
\be
\begin{split} \nonumber
AMSE(\hat{f},f) & \leq C_1(\gamma) \left\{\sum_{j \in {\cal J}_0^{c*}}
\left\{\sum_{|k| = k_j+1}^{(n_j-1)/2} |c_{kj}|^2+Pen_j(k_j)\right\}
+\sum_{j \in {\cal J}_0^*} Pen_j(1) + Pen_0(d)\right\} \\
& +C_2(\gamma) \frac{\sigma_2}{N} \{1-\pi_0(0)\}. \label{eq:case2}
\end{split}
\ee
We already showed that the first term $\sum_{j \in {\cal J}_0^{c*}}$
in the RHS of (\ref{eq:case2}) is
$O\left(\sum_{j \in {\cal J}_0^{c*}}\min(N^{-\frac{2s_j}{2s_j+1}},\frac{n_j}{N}\right)$.
The conditions of $\pi_j(1)$ and $\pi_0(d)$ imply that both
$\sum_{j \in J_0^*} Pen_j(1)$ and $Pen_0(d)$ are $O(d/N)$, and, therefore,
the first term in (\ref{eq:case2}) is dominating when $d^*_0 \sim d$.
\qed

\subsection*{Proof of Proposition \ref{prop:lower}} \label{subsec:proofproplow}
Consider the model (\ref{eq:barj}) and
the equivalent Gaussian sequence model (\ref{eq:fourier}) in the Fourier domain.
Evidently, $\inf_{\tilde{f}_j}\sup_{f_j \in \Fj}AMSE(\tilde{f}_j,f_j)=
\inf_{\tilde{c}_j}\sup_{c_j \in \Tj}E||\tilde{c}_j-c_j||^2_{n_j}$,
where $\tilde{c}_j$ are discrete Fourier coefficients of $\tilde{f}_j$.

Most of the proof is a direct consequence of the standard techniques for establishing
minimax lower bounds in the Gaussian sequence model over Sobolev ellipsoids
(see, e.g. Tsybakov, 2009, Section 3.2) but unlike the standard setup,
the variance in the considered model (\ref{eq:fourier}) depends on the sample
size $N$ that may affect the minimax rates.

Consider the class of diagonal linear estimators
$\tilde{c}_j(\lambda)$ of the form $\tilde{c}_{kj}=\lambda_k
\xi_{kj},\;k=-(n-1)_j/2,\ldots,-1,1,\ldots,(n_j-1)/2$ and
$\tilde{c}_{0j}=0$ (see Section \ref{subsec:idea}). It is well known
(see, e.g., Tsybakov, 2009, Section 3.2), that as $n_j$ tends to
infinity, the minimax linear diagonal estimator is asymptotically
minimax over all estimators of $f_j$: \be \nonumber
\begin{split}
\inf_{\tilde{c}_j}\sup_{c_j \in \Tj}
E||\tilde{c}_j-c_j||^2_{n_j} & \sim \inf_\lambda \sup_{c_j \in \Tj}
E||\tilde{c}_j(\lambda)-c_j||^2_{n_j} \\ & =
\sup_{c_j \in \Tj} \inf_\lambda E||\tilde{c}_j(\lambda)-c_j||^2_{n_j}.
\end{split}
\ee
By standard calculus (see, e.g., Tsybakov, 2009, Section 3.2),
\be \label{eq:inflam}
\inf_\lambda E||\tilde{c}_j(\lambda)-c_j||^2_{n_j}=\frac{\sigma^2}{N}
\sum_{k=-(n_j-1)/2}^{(n_j-1)/2}\frac{|c_{kj}|^2}{|c_{kj}|^2+\frac{\sigma^2}{N}}
\ee
and the minimax linear estimator $\hat{c}_j^L$ is then of the form
$$
\hat{c}^L_{kj}=(1-k^{s_j}\kappa_j)_+ \xi_{kj},
$$
where $\kappa_j$ is the solution of the equation
$$
\frac{\sigma^2}{N}\sum_{k=1}^{(n_j-1)/2}(2k)^{s_j}(1-(2k)^{s_j} \kappa_j)_+=\kappa_j R_j^2.
$$
Consider two cases:
\medskip
\noindent
\newline {\em a)} $2s_j+1 \geq \ln N/\ln n_j$. In this case we can follow Tsybakov
(2009, Section 3.2) to get
\be \label{eq:truncsum}
\frac{\sigma^2}{N}\sum_{k=1}^{(n_j-1)/2}(2k)^{s_j}(1-(2k)^{s_j} \kappa_j)_+=
\frac{\sigma^2}{N}\sum_{k=1}^{k_j}(2k)^{s_j}(1-(2k)^{s_j} \kappa_j),
\ee
where $k_j=\lfloor\frac{1}{2}\kappa_j^{-1/s_j}\rfloor$, and neglecting the constants,
$\kappa^2_j=N^{-\frac{2s_j}{2s_j+1}}$
and
$E||\hat{c}_j^L-c_j||^2_{n_j}=O\left(N^{-\frac{2s_j}{2s_j+1}}\right)$.

The condition $2s_j+1 \geq \ln N/\ln n_j$ is necessary to ensure
that the resulting $k_j= \frac{1}{2}N^{\frac{1}{2s_j+1}} \leq (n_j-1)/2$ in (\ref{eq:truncsum}).

\medskip
\noindent
\newline {\em b)} $2s_j+1<\ln N/\ln n_j$. In this case one can easily see that
$$
\frac{\sigma^2}{N}\sum_{k=1}^{(n_j-1)/2}(2k)^{s_j}(1-(2k)^{s_j} \kappa_j)_+=
\frac{\sigma^2}{N}\sum_{k=1}^{(n_j-1)/2}(2k)^{s_j}(1-(2k)^{s_j} \kappa_j),
$$
$\kappa_j^2=\frac{n_j}{N}$ and
$E||\hat{c}_j^L-c_j||^2_{n_j}=O\left(\frac{n_j}{N} \right)$.

\qed

\subsection*{Proof of Theorem \ref{th:lower}} \label{subsec:proofthlow}
No estimator $\tilde{f}$ of $f$ in (\ref{eq:spad}) can obviously perform
better than that of an oracle that knows the true subsets $\Jo$ and $\Jc$ of
zero and nonzero components $f_j$ of $f$. In this ideal case, one would certainly set $\hat{f}_j=0$
for all $j \in \Jo$ with no error and, therefore, due to the additivity of
the AMSE, Proposition \ref{prop:lower} yields
\be \nonumber
\begin{split}
\inf_{\tilde{f}}\sup_{f_j \in \Fj, j \in \Jc} AMSE(\tilde{f},f)
&=
\sum_{j \in \Jc} \inf_{\tilde{f}_j}\sup_{f_j \in \Fj} AMSE(\tilde{f}_j,f_j) \\
& \geq C_2 \sum_{j \in \Jc} \min\left(N^{-\frac{2s_j}{2s_j+1}},\frac{n_j}{N}\right)
\end{split}
\ee
(see Proposition 4.16 of Johnstone, 2013).

Furthermore, since $\min\left(N^{-\frac{2s_j}{2s_j+1}},\frac{n_j}{N}\right)
>N^{-1},\;j\in \Jc$, for $d_0 > d/2$ one has
$$
\frac{d_0 \ln(d/d_0)}{N} \leq \ln 2~ \frac{d_0}{N} \leq \ln 2 \sum_{j \in \Jc}
\min\left(N^{-\frac{2s_j}{2s_j+1}},\frac{n_j}{N}\right)
$$
and the first term in the RHS of (\ref{eq:lower}) is dominating.
Thus, to complete the proof we need to show that for $d_0 \leq d/2$,
\be \label{eq:kl} \inf_{\tilde{f}}\sup_{f_j \in \Fj, j \in \Jc}
AMSE(\tilde{f},f) = \inf_{\tilde{c}}\sup_{c_j \in \Tj,\;j \in \Jc}
||\tilde{c}-c||^2_{N_0} \geq C_2 \frac{d_0 \ln(d/d_0)}{N}, \ee where
$N_0=\sum_{j=1}^d n_j$ and $c$ is an $N_0$-dimensional amalgam of
$d$ $n_j$-dimensional vectors of discrete Fourier coefficients $c_j$
of $f_j$.

The proof is based on finding a subset ${\cal C}_{d_0}$ of $N_0$-dimensional
amalgamated vectors
$c$ with $d_0$ nonzero components $c_j \in \Tj$ such that for any pair
$c^1, c^2 \in {\cal C}_{d_0}$ and some constant $C>0$,
$||c^1-c^2||^2_{N_0} \geq C \frac{\sigma^2}{N}d_0 \ln(d/d_0)$ and the
Kullback-Leibler divergence $K(\mathbb{P}_{c^1},\mathbb{P}_{c^2})=
\frac{||c^1-c^2||^2_{N_0}}{2\sigma^2/N} \leq (1/16)\ln{\rm card}({\cal C}_{d_0})$.
The required result in (\ref{eq:kl}) then follows immediately from Lemma A.1
of Bunea et al. (2007).

Define the subset ${\cal \tilde{V}}_{d_0}$ of all $d$-dimensional indicator
vectors with $d_0$ entries of ones:
${\cal \tilde{D}}_{d_0}=\{v: v \in \{0,1\}^d,\; ||v||_0=d_0\}$.
Lemma A.3 of Rigollet \& Tsybakov (2011) implies that for $d_0 \leq d/2$,
there exists a subset
${\cal V}_{d_0} \subset {\cal \tilde{V}}_{d_0}$ such that for some constant
$C_0>0$, $\ln {\rm card}({\cal V}_{d_0}) \geq C_0 d_0\ln(d/d_0)$,
and for any pair $v_1,v_2 \in {\cal V}_{d_0}$, the Hamming distance
$\rho(v_1,v_2)=\sum_{j=1}^d \mathbb{I}\{v_{1j} \neq v_{2j}\} \geq C_0
d_0$.

To any indicator vector $v \in {\cal V}_{d_0}$ assign the
corresponding vector $c \in {\cal C}_{d_0}$  as follows. Let
$\tilde{C}^2=(1/16)C_0 \frac{\sigma^2}{N} \ln(d/d_0)$. Define $c_j$ to be a
zero
vector if $v_j=0$ and to have two nonzero entries $c_{-1j}=c_{1j}=\tilde{C}/\sqrt{2}$ otherwise.
Evidently, $c_j \in {\cal F}(s_j, \tilde{C}) \subset \Fj$ for $v_j=1$
and ${\rm card}({\cal C}_{d_0})={\rm card}({\cal V}_{d_0})$.

For any pair $c^1, c^2 \in {\cal C}_{d_0}$ and the corresponding
$v_1, v_2 \in {\cal V}_{d_0}$, we then have
\begin{equation*}
\begin{split}
||c^1-c^2||^2_{N_0}& = \tilde{C}^2 \sum_{j=1}^d
\mathbb{I}\{v_{1j} \neq v_{2j}\} \geq \tilde{C}^2 \; C_0\;
d_0 = \frac{1}{16}\frac{\sigma^2}{N} C_0^2 d_0\ln(d/d_0), \\
K(\mathbb{P}_{c^1},\mathbb{P}_{c^2}) &=
\frac{\tilde{C}^2}{2\sigma^2/N} \sum_{j=1}^d
\mathbb{I}\{v_{1j} \neq v_{2j}\} \leq \frac{\tilde{C}^2
d_0}{\sigma^2/N} \leq \frac{1}{16}\ln{\rm card}({\cal C}_{d_0}),
\end{split}
\end{equation*}
which completes the proof.
\qed

\end{document}